\documentclass[12pt,a4paper]{amsart}

\usepackage[utf8]{inputenc}
\usepackage[margin=1in]{geometry}
\usepackage{thmtools} %Manipulating ref package
\usepackage{hyperref}
\usepackage{etoolbox} %Clever Refencing inside enumerate environment.
\usepackage[capitalize,noabbrev]{cleveref}
\usepackage{amsmath}
\usepackage{amsfonts}
\usepackage{amsthm}
\usepackage{amssymb}
\usepackage{xcolor}
\usepackage{graphicx}
\usepackage{tikz-cd} %Diagrams
\usepackage{bbm} %Little case k
\usepackage{enumerate} %Enumerate environment
\usepackage{stmaryrd} %Mapsfrom
\usepackage{subcaption} %Subfigures
\usepackage{microtype} %Changes appearance font
\usepackage{mathrsfs} %Fonts e.g. mathscr
\usepackage{enumitem}
\usepackage{comment}

\theoremstyle{plain}
\newtheorem{thrm}{Theorem}[section]
\def\bthm{\begin{thrm}}
\def\ethm{\end{thrm}}

\newtheorem{prop}[thrm]{Proposition}
\def\bprop{\begin{prop}}
\def\eprop{\end{prop}}

\newtheorem{ques}[thrm]{Question}
\def\bques{\begin{ques}}
\def\eques{\end{ques}}

\newtheorem{cjec}[thrm]{Conjecture}
\def\bcjec{\begin{cjec}}
\def\ecjec{\end{cjec}}

\newtheorem{cor}[thrm]{Corollary}
\def\bcor{\begin{cor}}
\def\ecor{\end{cor}}

\newtheorem{fact}[thrm]{Fact}
\def\bfact{\begin{fact}}
\def\efact{\end{fact}}

\newtheorem{lem}[thrm]{Lemma}
\def\blem{\begin{lem}}
\def\elem{\end{lem}}

\theoremstyle{definition}
\newtheorem{defn}[thrm]{Definition}
\def\bdefn{\begin{defn}}
\def\edefn{\end{defn}}

\newtheorem{nota}[thrm]{Notation}
\def\bnota{\begin{nota}}
\def\enota{\end{nota}}

\newtheorem*{conc}{Conclusion}
\def\bconc{\begin{conc}}
\def\econc{\end{conc}}

\newtheorem{alg}[thrm]{Algorithm}
\def\balg{\begin{alg}}
\def\ealg{\end{alg}}

\def\bproof{\begin{proof}}
\def\eproof{\end{proof}}

\theoremstyle{remark}

\newtheorem{rem}[thrm]{Remark}
\def\brem{\begin{rem}}
\def\erem{\end{rem}}

\newtheorem{ex}[thrm]{Example}
\def\bex{\begin{ex}}
\def\eex{\end{ex}}

\newtheorem{exs}[thrm]{Examples}
\def\bexs{\begin{exs}}
\def\eexs{\end{exs}}

\newtheorem{obs}{Observation}
\def\bobs{\begin{obs}}
\def\eobs{\end{obs}}	
\let\ra=\rightarrow

\let\setminus=\smallsetminus

\let\epsilon=\varepsilon

\def\A{\mathcal{A}}
\def\B{\mathcal{B}}
\def\H{\mathcal{H}}

\def\k{\mathbbm{k}}

\DeclareMathOperator{\hgt}{ht}

\DeclareMathOperator{\Tor}{Tor}

\DeclareMathOperator{\init}{in}

\DeclareMathOperator{\reg}{reg}

\DeclareMathOperator{\im}{im}

\DeclareMathOperator{\iv}{iv}
\DeclareMathOperator{\St}{St}

\newcommand{\ul}[1]{\underline{#1}}

\renewcommand{\subset}{\subseteq}

\AtBeginDocument{%
   \def\MR#1{}
}

\begin{document}

\title[Castelnuovo--Mumford Regularity]{A New Conjecture and Upper Bound\\ on the Castelnuovo--Mumford Regularity\\ of Binomial Edge Ideals}
\author{Adam LaClair}
%\thanks{The author was partially supported by NSF grant DMS-2100288 and by Simons Foundation Collaboration Grant for Mathematicians \#580839}

%    \subjclass is required.

\subjclass[2020]{Primary 13D02, 13F70. Secondary 05E40, 13D22.}

\begin{abstract}
A famous theorem of Kalai and Meshulam is that $\reg(I + J) \leq \reg(I) + \reg(J) -1$ for any squarefree monomial ideals $I$ and $J$. This result was subsequently extended by Herzog to the case where $I$ and $J$ are any monomial ideals. In this paper we conjecture that the Castelnuovo--Mumford regularity is subadditive on binomial edge ideals. Specifically, we propose that $\reg(J_{G}) \leq \reg(J_{H_{1}}) + \reg(J_{H_{2}}) -1$ whenever $G$, $H_{1}$, and $H_{2}$ are graphs satisfying $E(G) = E(H_{1}) \cup E(H_{2})$ and $J_{\ast}$ is the associated binomial edge ideal. We prove a special case of this conjecture which strengthens the celebrated theorem of Malayeri--Madani--Kiani that $\reg(J_{G})$ is bounded above by the minimal number of maximal cliques covering the edges of the graph $G$. From this special case we obtain a new upper bound for $\reg(J_{G})$, namely that $\reg(J_{G}) \leq \hgt(J_{G}) +1$. Our upper bound gives an analogue of the well-known result that $\reg(I(G)) \leq \hgt(I(G)) +1$ where $I(G)$ is the edge ideal of the graph $G$. We additionally prove that this conjecture holds for graphs admitting a combinatorial description for its Castelnuovo--Mumford regularity, that is for closed graphs, bipartite graphs with $J_{G}$ Cohen--Macaulay, and block graphs. Finally, we give examples to show that our new upper bound is incomparable with Malayeri--Madani--Kiani's upper bound for $\reg(J_{G})$ given by the size of a maximal clique disjoint set of edges.
\end{abstract}

\keywords{Castelnuovo--Mumford regularity, binomial edge ideals, Terai's conjecture, Eisenbud--Goto conjecture.}

\dedicatory{Dedicated to the memory of Professor J\"urgen Herzog}

\maketitle

\section{Introduction}

Castelnuovo--Mumford regularity introduced by David Mumford in \cite{mumford1966lectures} is an important invariant appearing in commutative algebra and algebraic geometry which gives bounds on the non-vanishing of the graded components of $\Tor$ modules and local cohomology modules supported at the maximal ideal of a graded module over a polynomial ring. Unlike projective dimension which is bounded by the Krull dimension of the polynomial ring by Hilbert's syzygy theorem, there are examples \cite{mayr1982complexity}, \cite{Koh1998Quadrics} which show that Castelnuovo--Mumford regularity of an ideal can be doubly exponential in the Krull dimension of the polynomial ring even when all of the generators have low degree and the number of minimal generators of the ideal is \textit{small}. The Eisenbud--Goto conjecture \cite{Eisenbud1984Regularity} famously proposed a linear bound on the Castelnuovo--Mumford regularity of non-degenerate prime ideals in terms of the degree of the ideal which has since been shown not to hold by counterexamples of McCullough and Peeva \cite{McCullough2018Counterexamples}, see also \cite{Caviglia2019RegularityPrime}. The question of bounding the Castelnuovo--Mumford regularity of ideals in terms of other homological invariants continues to draw much attention. 

The starting point for this paper is the famous conjecture of Terai \cite{terai2001EisenbudGoto}, arising from his study of the Eisenbud--Goto conjecture, which asserted that Castelnuovo--Mumford regularity is subadditive along squarefree monomial ideals. This conjecture was proved by Kalai and Meshulam \cite{kalai2006Intersections}, and subsequently generalized to arbitrary monomial ideals by Herzog \cite{herzog2007Generalization}. Specifically, it has been shown that
\bthm[\cite{herzog2007Generalization}, \cite{kalai2006Intersections}, \cite{terai2001EisenbudGoto}]
\label{thm:kalai}
Let $I$ and $J$ be monomial ideals of a polynomial ring $S = \k[x_{1},\ldots,x_{n}]$, then 
\begin{align*}
    \reg(S/(I+J)) \leq \reg(S/I) + \reg(S/J).
\end{align*}
\ethm
A consequence of Theorem \ref{thm:kalai} is that if $G$ is a graph and $I(G)$ is the associated edge ideal, then $\reg(I(G)) \leq \hgt(I(G))+1$. For edge ideals this upper bound can be improved from $\hgt(I(G))$, which is equal to the size of a minimal vertex cover, to the size of a maximal matching, see \cite{Ha2008Monomial}. Inspired by these results, we investigate whether an analogue of Theorem \ref{thm:kalai} holds for binomial edge ideals. Based on our work in this paper we are led to conjecture the following:
\bcjec
\label{conj:subadditivity}
Let $G$ be a graph on $n$ vertices, let $R := \k[x_{1},\ldots,x_{n},y_{1},\ldots,y_{n}]$ for any field $\k$, and define the binomial edge ideal $$J_{G} := ( x_{i}y_{j} - x_{j}y_{i} \mid \{i,j\} \in E(G) ).$$ If $H_{1}$ and $H_{2}$ are subgraphs of $G$ satisfying $E(G) = E(H_{1}) \cup E(H_{2})$, then
\begin{align*}
    \reg(R/J_{G}) \leq \reg(R/J_{H_{1}}) + \reg(R/J_{H_{2}}).
\end{align*}
\ecjec

A special case of Conjecture \ref{conj:subadditivity} is given by the theorem of Maleyeri--Madani--Kiani which bounds the regularity of the binomial edge ideal from above by the number of maximal cliques of the graph. 
\bthm[\cite{malayeri2021proof}]
\label{thm:reg_upper_bnd_malayeri}
Let $G$ be a graph, and let $c(G)$ denote the number of maximal cliques of $G$. In particular, $E(G)$ can be covered by $c(G)$ complete graphs. Then, 
\begin{align*}
    \reg(R/J_{G}) \leq c(G).
\end{align*}
\ethm 
In fact, the authors gave a sharper upper bound for $\reg(J_{G})$ in terms of the size of a maximal clique disjoint set of edges \cite{malayeri2021proof}. 

The first result of this paper is a generalization of Theorem \ref{thm:reg_upper_bnd_malayeri}. We prove
\bthm
\label{thm:subadditivity_special_case}
Let $G$ be a graph. Let $\A$ and $\B$ be a set of subgraphs of $G$ such that every element of $\A$ is a complete graph of the form $K_{m_{i}}$ for some $m_{i} \geq 2$, and every element of $\B$ is a complete bipartite graph of the form $K_{1,n_{i}}$ for some $n_{i} \geq 2$. Suppose that $\A\cup \B$ is an edge cover of $G$, then
\begin{align*}
    \reg(R/J_{G}) \leq \#\A + 2 \cdot \#\B.
\end{align*}
\ethm
From Theorem \ref{thm:subadditivity_special_case} we deduce the following result.
\bthm
\label{thm:reg_leq_ht}
Let $G$ be a graph, then 
\begin{align*}
    \reg(R/J_{G}) \leq \hgt(J_{G}).
\end{align*}
\ethm 
Theorem \ref{thm:reg_leq_ht} strengthens a result of Matsuda--Murai \cite{matsuda2013regularity} that $\reg(R/J_{G}) \leq n-1$ where $n$ is the number of vertices of $G$. Indeed, letting $K_{n}$ denote the complete graph on $n$ vertices we have that $J_{K_{n}}$ is a minimal prime of $J_{G}$ \cite{herzog2010binomial}, and hence $\hgt(J_{G}) \leq \hgt(J_{K_{n}}) = n-1$. Our proof of Theorems \ref{thm:subadditivity_special_case} and \ref{thm:reg_leq_ht} go via application of the regularity lemma, Lemma \ref{lem:regularity_lemma}, applied to short exact sequences which relate the Castelnuovo--Mumford regularity of $J_{G}$ to the Castelnuovo--Mumford regularity of auxillary graphs appearing in a short exact sequence. We review these lemmas in the next section.

Lastly, we show that Conjecture \ref{conj:subadditivity} holds for closed graphs, bipartite graphs having Cohen--Macaulay binomial edge ideals, and block graphs. This result follows from previous work of the author and others (\cite{ene2015regularity}, \cite{jayanthan2019regularityCM&BipartiteGraphs}, \cite{laclair2024regularity}) establishing a combinatorial interpretation of Castelnuovo--Mumford regularity for these graphs in terms of lengths of certain induced subpaths of the graph.

\section{Background}

For us by (simple) graph $G$ we mean a pair of vertices $V$ and edges $E$ such that $E$ contains no loops nor multi-edges. A set $\A$ consisting of subgraphs of $G$ is called an \textbf{edge cover of } $\mathbf{G}$ if $E(G) = \bigcup_{H \in \A} E(H)$. For a vertex $v \in G$ we define the \textbf{neighborhood of $\mathbf{v}$} as $$N_{G}(v) := \{w \in G \mid \{v,w\} \in E(G)\}.$$ We define the \textbf{closed neighborhood of $\mathbf{v}$} as $$N_{G}[v] := N_{G}(v) \cup \{v\}.$$ For terminology on graph theory not defined we refer the reader to \cite{west1996graphTheory}. The complete bipartite graph of the form $K_{1,m}$ for $m \geq 1$ is defined to have vertex set $\{1,2,\ldots,m,m+1\}$ and edge set $\{ \{1,i\} \mid 2 \leq i \leq m+1\}$. We will call the vertex $1$ to be the \textbf{center} of $K_{1,m}$. Given a vertex $v \in G$ with at least two neighbors we denote the \textbf{star} of $G$ at $v$, denoted $\St_{G}(v)$, to be the subgraph of $G$ with $$V(\St_{G}(v)) := N_{G}[v]$$ and $$E(\St_{G}(v)) := \{ \{a,v\} \mid a \in N_{G}(v) \}.$$

\bdefn
For a graph $G$ say that a vertex $v\in G$ is \textbf{simplicial} or \textbf{free} if the induced subgraph on $N_{G}[v]$ is a complete graph. For a graph $G$ denote the number of non-simplicial vertices of $G$ by $\iv(G)$. Observe that $\iv(G) = \sum_{i=1}^{c(G)} \iv(G_{i})$ whenever $G = \bigsqcup_{i=1}^{c(G)} G_{i}$. % For a subset $A$ of $G$ we denote by $\iv(A)$ the number of vertices of $A$ which are not simplicial vertices of $G$.
\edefn 

Given a graph $G$ and vertex $v$ we denote by
\begin{itemize}
    \item $G\setminus v$ the subgraph having vertices $V(G) \setminus v$ and edge set $E(G) \setminus \{ \{u,v\} \mid u \in N_{G}(v) \}$, and
    \item $G_{v}$ the graph having vertices $V(G)$ and edge set $E(G) \cup \{ \{a,b\} \mid a,b\in N_{G}(v)\}$.
\end{itemize}

Recall the following fundamental result of Ohtani.
\blem[{\cite[Lemma 4.8]{ohtani2011graphs}}]
\label{lem:exact_sequence_binomial_edge_ideal}
Let $G$ be a graph and $v \in G$ a non-simplicial vertex. Then, there is a short exact sequence of $R$-modules
\begin{equation}
\label{eqn:ses_binomial_edge_ideal}
    0 \ra R/J_{G} \ra R/J_{G_{v}} \oplus R_{v}/J_{G \setminus v} \ra R_{v}/J_{G_{v} \setminus v} \ra 0
\end{equation}
where $R_{v} = \k[X_{i},Y_{i} : i \in [n] \setminus \{v\}]$.
\elem

Recall the well-known regularity lemma for a short exact sequence.
\blem[cf. {\cite[Corollary 18.7]{peeva2010graded}}]
\label{lem:regularity_lemma}
Suppose that $0 \ra U \ra U^{'} \ra U^{''} \ra 0$ is a short exact sequence of graded finitely generated $R$-modules with homomorphisms of degree $0$. Then,
\begin{enumerate}
    \item $\reg(U^{'}) \leq \max\{ \reg(U), \reg(U^{''})\}$,
    \item $\reg(U) \leq \max\{ \reg(U^{'}), \reg(U^{''})+1\}$
    \item $\reg(U^{''}) \leq \max\{ \reg(U)-1, \reg(U^{'})\}$
\end{enumerate}
\elem 

%Lemma that I don't believe that I need.
\begin{comment}
\blem 
Let $G$ be a graph, $v$ a simplicial vertex of $G$, and $w$ any vertex of $G$. Then, $v$ is a simplicial vertex of $G_{w}$.
\elem 

\bproof
We may suppose that $v \in N_{G}(w)$; otherwise, there is nothing to prove. In which case $N_{G_{v}}[w] = N_{G}[w] \cup N_{G}[v]$. By definition of $G_{v}$ and the fact that $N_{G}[w]$ is complete, it follows that $N_{G_{v}}[w]$ is complete.
\eproof
\end{comment}

\section{Results}

We are now ready to prove Theorem \ref{thm:subadditivity_special_case}.

\bthm
\label{thm:covering_result}
Let $G$ be a graph. Let $\A := \{A_{1},\ldots,A_{p}\}$ and $\B := \{B_{1},\ldots,B_{q}\}$ be sets of subgraphs of $G$ such that for all $1 \leq i \leq p$, $A_{i} = K_{n_{i}}$ for some $n_{i} \geq 2$, and for all $1 \leq i \leq q$, $B_{i} = K_{1,m_{i}}$ for some $m_{i} \geq 2$. Suppose that $\A \cup \B$ is an edge cover of $G$. Then,
\begin{align*}
\reg(R/J_{G}) \leq p + 2 \cdot q.
\end{align*}
\ethm 

\bproof
We induce on $q$. If $q = 0$, then $c(G) \leq p$, and the theorem statement follows from Theorem \ref{thm:reg_upper_bnd_malayeri}. 

Suppose that $q > 0$, and suppose by induction hypothesis that for every graph $G^{'}$ admitting an edge covering by $(\A^{'},\B^{'})$ with $\A^{'}$ and $\B^{'}$ defined analogously as in the theorem statement, and with $\# \B^{'} < q$ that $\reg(R/J_{G^{'}}) \leq \# \A^{'} + 2 \# \B^{'}$. Let $G$ be a graph admitting an edge covering by $(\A,\B)$ as in the theorem statement with $\# \B = q$. Denote by $v$ the center of $B_{q}$. Observe that if $v$ is a simplicial vertex of $G$, then the pair $\A^{'} := \A \cup \{ N_{G}[v] \}$ and $\B^{'} := \B \setminus \{B_{q}\}$ would realize an edge cover of $G$ with $\# \B^{'} < \# \B$ and satisfy $\# \A^{'} + 2 \#\B^{'} < \# \A + 2 \#\B$. Hence the theorem statement would follow from the induction hypothesis. Thus, we may assume that $v$ is not a simplicial vertex of $G$. We may further replace $B_{q}$ by the star of $G$ at $v$, as this will not change the value $\# \A + 2\#\B$. Obviously, we may assume that $v$ is not the center of $B_{i}$ for $1 \leq i \leq q-1$. Otherwise, we could replace $\B$ by $\B \setminus \{B_{i}\}$ and conclude by induction hypothesis since $B_{i}$ would be a subgraph of $B_{q}$. We may assume that $\# (B_{i} \setminus v) \neq 2$ for $1 \leq i \leq q-1$. Otherwise, $B_{i} \setminus v$ is a complete graph, and $\A^{'} := \A \cup \{B_{i} \setminus v \}$, $\B^{'} := \B \setminus B_{i}$ is an edge cover of $G$ with $\# \B^{'} < \# \B$ and $\# \A^{'} + 2\# \B^{'} < \# \A + 2\# \B$, and we could conclude by induction hypothesis. Similarly, we may assume that $\# (A_{i} \setminus v) \geq 2$ for $1 \leq i \leq p$. In summary, we have reduced to the following setup:
\begin{itemize}
    \item $v := v_{q}$ is not not a simplicial vertex of $G$,
    \item $B_{i} \setminus v$ is a star of $G$ on at least three vertices for $1 \leq i \leq q-1$,
    \item $A_{i} \setminus v$ is a complete subgraph of $G$ on at least two vertices for $1 \leq i \leq p$.
\end{itemize}
After the above reductions we may observe that
\begin{itemize}
    \item $G \setminus v$ has an edge cover given by $\A^{'} := \{ A_{i} \setminus v \mid 1 \leq i \leq p \}$ and $\B^{'} := \{ B_{i} \setminus v \mid 1 \leq i \leq q-1 \}$,
    \item $G_{v}$ has an edge cover given by $\A^{'} := \A \cup \{ (B_{q})_{v} \} $ and $\B^{'} := \B \setminus \{B_{q}\}$,
    \item $G_{v} \setminus v$ has an edge cover given by $\A^{'} := \{ A_{i} \setminus v \mid 1 \leq i \leq p \} \cup \{ (B_{q})_{v} \setminus v \}$ and $\B^{'} := \{ B_{i} \setminus v \mid 1 \leq i \leq q-1 \}$. (Observe that $\# (B_{q})_{v} \setminus v \geq 2$ since $\# B_{q} = m_{q} +1 \geq 3$.)
\end{itemize}
In the above instances we have an edge covering of a graph with smaller $\# \B^{'} < \# \B = q$. The induction hypothesis together with Lemmas \ref{lem:exact_sequence_binomial_edge_ideal} and \ref{lem:regularity_lemma} imply that
\begin{align*}
    \reg(R/J_{G}) \leq \max \{ p + 2(q-1), p+1 + 2(q-1), \big( p+ 1 + 2(q-1) \big) + 1\} = p + 2q.
\end{align*}
\eproof

Recall the following description for the minimal primes of $J_{G}$.

\bprop[{\cite[Lemma 3.1, Theorem 3.2, Corollary 3.9]{herzog2010binomial}}]
\label{prop:min_primes_binomial_edge_ideal}
Let $G$ be a connected simple graph on $[n]$. For a subset $S \subset [n]$, define 
\begin{align*}
    P_{G}(S) := (\bigcup_{i \in S} \{x_{i},y_{i}\},J_{\tilde{G}_{1}},\ldots,J_{\tilde{G}_{c(S)}})
\end{align*}
where $G_{1},\ldots,G_{c(S)}$ denote the connected components of $G \setminus S$ and $\tilde{G}_{i}$ denotes the complete graph on the vertex set $V(G_{i})$. Then,
\begin{enumerate}
    \item $P_{G}(S)$ is a prime ideal,
    \item $\hgt P_{G}(S) = 2 \#S + (n - c(S))$,
    \item $J_{G} = \bigcap_{S \subset [n]} P_{G}(S)$,
    \item $P_{G}(S)$ is a minimal prime of $J_{G}$ if and only if $S = \varnothing$, or $S \neq \varnothing$ and for each $i \in S$ one has $c(S\setminus \{i\}) < c(S)$. \label{item:min_primes_3}
\end{enumerate}
\eprop

\brem
We say that $S$ satisfying Item \ref{item:min_primes_3} is a cut set of $G$.
\erem

We introduce the following notation. Let $S \subset V(G)$ (not necessarily a cut set of $G$), and $G_{1},\ldots,G_{c}$ denote the connected components of $G \setminus S$. Define $$b_{G}(S) := 2 \# S + \sum_{i=1}^{c} (\# G_{i} - 1).$$ Define the invariant 
\begin{align*}
b_{G} := \min\{b_{G}(S) : S \subset V(G) \}.
\end{align*}

\blem
\label{lem:b_G_equals_height}
For a graph $G$ we have that
\begin{align*}
b_{G} = \hgt J_{G}.
\end{align*}
\elem

\bproof
When $S$ is a cut set of $G$ it follows from the proof of {\cite[Lemma 3.1]{herzog2010binomial}} that $b_{G}(S) = \hgt P_{G}(S)$. Choosing $S$ such that $\hgt P_{G}(S) = \hgt J_{G}$, we have that 
\begin{align*}
    b_{G} \leq b_{G}(S) = \hgt P_{G}(S) = \hgt J_{G}.
\end{align*}

Next, we show that $\hgt J_{G} \leq b_{G}$. It suffices to show that if $S \subset V(G)$ is chosen such that $b_{G} = b_{G}(S)$, then $S$ is a cut set of $G$. Suppose by contradiction that $S$ is not a cut set of $G$. Then, $S \neq \varnothing$, and there exists $i \in S$ such that $c(G \setminus (S \setminus \{i\})) = c(G \setminus S)$. If $G_{1},\ldots,G_{c}$ denote the connected components of $G \setminus S$, then $i$ is adjacent to exactly one of the connected components, say $G_{1}$. Put $\tilde{S} = S \setminus \{i\}$. Then, $G \setminus \tilde{S}$ has connected components $G_{1} \cup \{i\}, G_{2},\ldots, G_{c}$. It follows that 
\begin{align*}
b_{G}(\tilde{S}) < b_{G}(S)
\end{align*}
which contradicts our choice of $S$. Thus, $S$ must be a cut set of $G$.
\eproof

We are now ready to prove the main result of this paper.

\bthm
For a graph $G$ we have that
\begin{align*}
    \reg(R/J_{G}) \leq \hgt(J_{G}).
\end{align*}
\ethm

\bproof
We may assume that $G$ is connected since Castelnuovo--Mumford regularity and height is additive over connected components. We may assume that $G$ is not a vertex; otherwise, $J_{G} = 0$. We induce on $\iv(G)$. If $\iv(G) = 0$, then $G$ is a complete graph. In which case, $J_{G}$ is equal to the ideal of two minors of a generic $2 \times n$ matrix, and $\reg(R/J_{G}) = 1$ and $\hgt J_{G} \geq 1$. Suppose that $\iv(G) > 0$ and that the claim has been shown for all graphs having strictly smaller number of non-simplicial vertices.

Pick $S \subset V(G)$ such that $b_{G} = b_{G}(S)$. Let $G_{1},\ldots,G_{c}$ denote the connected components of $G \setminus S$ which are not isolated vertices.

\ul{Case 1.} Suppose that $G_{i}$ is a complete graph for all $1 \leq i \leq c$.

\textit{Proof of Case 1.} Put $\A := \{G_{1},\ldots,G_{c}\}$, and $\B := \{ \St_{G}(v) \mid v \in S \}$. Then $\A \cup \B$ realizes an edge cover of $G$. By Theorem \ref{thm:covering_result} we have that 
\begin{align*}
    \reg(R/J_{G}) \leq c + 2 \# S \leq b_{G}(S) = b_{G}.
\end{align*}

\ul{Case 2.} Without loss of generality we may assume that $G_{1}$ is not a complete graph.

\textit{Proof of Case 2.} It follows that there exists a non-simplicial vertex $v$ of $G_{1}$, and in particular $v$ is not a simplicial vertex of $G$. Let $H_{1},\ldots,H_{d}$ denote the connected components of $G_{1} \setminus v$ for some $d \geq 1$. Then, the connected components of $(G\setminus v) \setminus S$ are $\{H_{1},\ldots,H_{d},G_{2},\ldots,G_{c}\}$. It follows that
\begin{align*}
    b_{G \setminus v}(S) &= 2\# S + \sum_{i=1}^{d} (\# H_{i} - 1) + \sum_{i=2}^{c} (\# G_{i} - 1) \\
    &= 2\# S + (\# G_{1} - 1) -d + \sum_{i=2}^{c} (\# G_{i} - 1) \\
    &< b_{G}(S).
\end{align*}
The connected components of $G_{v} \setminus S$ are $\{(G_{1})_{v},G_{2},\ldots,G_{c}\}$, and it follows that $b_{G_{v}}(S) = b_{G}(S)$.
The connected components of $(G_{v} \setminus v) \setminus S$ are $\{(G_{1})_{v} \setminus v,G_{2},\ldots,G_{c}\}$ because $(G_{1})_{v} \setminus v$ is connected. Hence, it follows that $b_{G_{v} \setminus v}(S) < b_{G}(S)$. 

Thus, we compute that
\begin{align*}
    \reg(R/J_{G}) 
    &\leq \max\{ \reg(R/J_{G\setminus v}),\reg(R/J_{G_{v}}), \reg(R/J_{G_{v}\setminus v}) +1\} \\
    &\leq \max\{ b_{G \setminus v}(S), b_{G_{v}}(S), b_{G_{v}\setminus v}(S) +1 \} \\
    &\leq b_{G}(S) \\
    &= b_{G} \\
    &= \hgt J_{G}.
\end{align*}
The first inequality is Lemmas \ref{lem:exact_sequence_binomial_edge_ideal}, \ref{lem:regularity_lemma}. The second inequality follows from induction hypothesis on the number of non-simplicial vertices together with {\cite[Lemma 2.2]{malayeri2021proof}} which is the statement that $\max\{\iv(G \setminus v), \iv(G_{v}), \iv(G_{v} \setminus v)\} < \iv(G)$. The third inequality follows from the computations in the preceeding paragraph. The fourth and fifth equality follow from our choice of $S$ and Lemma \ref{lem:b_G_equals_height}.
\eproof

\brem
Proposition \ref{prop:min_primes_binomial_edge_ideal} shows that the height of the binomial edge ideal is combinatorially determined. Thus, $\hgt J_{G}$ gives a combinatorial upper bound for $\reg(R/J_{G})$.
\erem

In \cite{malayeri2021proof} Malayeri--Madani--Kiani prove that for any graph $G$ that $\reg(R/J_{G}) \leq \eta(G)$ where $$\eta(G) := \max\{ \# \H \mid \H \text{ is a clique disjoint edge set of } G\}.$$ We show that $\eta(G)$ and $\hgt J_{G}$ are incomparable. In other words there exists graphs $G_{1}$ and $G_{2}$ such that $\eta(G_{1}) < \hgt J_{G_{1}}$ and $\hgt J_{G_{2}} < \eta(G_{2})$. When $G$ is the net, see Figure \ref{fig:the_net}, we have that $\reg(R/J_{G}) = \eta(G) = 4 < \hgt J_{G} = 5$. When $G = K_{1,m}$ for some $m \geq 3$, we have that $\reg(R/J_{G}) = \hgt J_{G} = 2 < \eta(G) = m$.

\def\putBelow{-.5*.9}
\def\edgeDist{1.75*.9}
\def\rtTriangleEdgeMultiplier{1/8*.9}

\begin{figure}[h]
\begin{center}
\begin{tikzpicture}
%Drawing of G~
\filldraw[black] (-3*\edgeDist,0) circle (2pt) node at (-3*\edgeDist,\putBelow) {$1$};
\filldraw[black] (-2*\edgeDist,0) circle (2pt) node at (-2*\edgeDist,\putBelow) {$2$};
\filldraw[black] (-1*\edgeDist,0) circle (2pt) node at (-1*\edgeDist,\putBelow) {$3$};
\filldraw[black] (0,0) circle (2pt) node at (0,\putBelow) {$4$};
%\filldraw[black] (1*\edgeDist,0) circle (2pt) node at (1*\edgeDist,\putBelow) {$7$};
%\filldraw[black] (2*\edgeDist,0) circle (2pt) node at (2*\edgeDist,\putBelow) {$8$};
%\filldraw[black] (3*\edgeDist,0) circle (2pt) node at (3*\edgeDist,\putBelow) {$9$};
\filldraw[black] (-3/2*\edgeDist, 
0.866025403784*\edgeDist) circle (2pt) node at (-3/2*\edgeDist+.5, 
0.866025403784*\edgeDist) {$5$};
\filldraw[black] (-3/2*\edgeDist,1.75*\edgeDist) circle (2pt) node at (-3/2*\edgeDist+.5,1.75*\edgeDist) {$6$};

\draw (-3*\edgeDist,0) to (-2*\edgeDist,0);
\draw (-2*\edgeDist,0) to (-1*\edgeDist,0);
\draw (-1*\edgeDist,0) to (0,0);
\draw (-2*\edgeDist,0) to (-3/2*\edgeDist, 
0.866025403784*\edgeDist);
\draw (-1*\edgeDist,0) to (-3/2*\edgeDist, 
0.866025403784*\edgeDist);
\draw (-3/2*\edgeDist, 
0.866025403784*\edgeDist) to (-3/2*\edgeDist,1.75*\edgeDist);

\end{tikzpicture}
\end{center}
\caption{Net}
\label{fig:the_net}
\end{figure}

Next, we show that Conjecture \ref{conj:subadditivity} holds for certain families of graphs admitting combinatorial descriptions for $\reg(R/J_{G})$.

\bthm
Conjecture \ref{conj:subadditivity} holds for $G$ belonging to any one of the following families of graphs:
\begin{enumerate}
\item $G$ is a closed graph,
\item $G$ is a bipartite graph with $J_{G}$ Cohen--Macaulay,
\item $G$ is a block graph
\end{enumerate} 
\ethm

\bproof
Let $H_{1}$ and $H_{2}$ realize an edge cover of $G$. By {\cite[Theorems 4.1, 5.17]{laclair2024regularity}} we have that for $G$ as specified above that $\nu(G) = \reg(R/J_{G})$, and moreover that $\nu(G) = \#E(P)$ where $P$ is a union of vertex disjoint paths of $G$. Let $P_{1}$ and $P_{2}$ denote the restriction of $P$ to $H_{1}$ and $H_{2}$, respectively. By {\cite[Definition 3.4, Theorem 3.15]{laclair2024regularity}} we have that $\# E(P_{i}) \leq \nu(H_{i})$ for $i = 1,2$. Consequently, we have that
\begin{align*}
\reg(R/J_{G}) = \nu(G) = \# E(P) \leq \# E(P_{1}) + \# E(P_{2}) \leq \nu(H_{1}) + \nu(H_{2}) \leq \reg(R/J_{H_{1}}) + \reg(R/J_{H_{2}})
\end{align*}
where the last inequality follows from {\cite[Proposition 3.9]{laclair2024regularity}}.
\eproof

Lastly, we prove a special case of Conjecture \ref{conj:subadditivity} holds when the intersection of the graphs $H_{1}$ and $H_{2}$ is particularly ``nice".

\bthm
\label{thm:graph_decomp_implies_conjec}
Let $G$ be a connected graph, and $H_{1}$ and $H_{2}$ subgraphs of $G$ such that $E(G) = E(H_{1}) \cup E(H_{2})$. Furthermore, suppose that there are non-empty subsets of $V(G)$, $A$, $B$, and $C$, satisfying the following:
\begin{enumerate}
    \item $V(G) = A \sqcup B \sqcup C$,
    \item the subgraph of $G$ induced by $B$ is a connected complete subgraph of $G$, \label{item:second_item_in_thm_graph_decomp}
    \item $H_{1}$ is equal to the subgraph of $G$ induced by $A \cup B$, and $H_{2}$ is the subgraph of $G$ induced by $B \cup C$. \label{item:third_item_in_thm_graph_decomp}
\end{enumerate}
Then,
\begin{align*}
    \reg(R/J_{G}) \leq \reg(R/J_{H_{1}}) + \reg(R/J_{H_{2}}).
\end{align*}
\ethm

\bproof
Without loss of generality we may suppose that $G$ is a connected graph having at least two edges. We may suppose that $H_{1} \not\subset H_{2}$ and that $H_{2} \not\subset H_{1}$. Label the vertices of $G$ such that the following conditions are satisfied:
\begin{itemize}
    \item $i < j$ if $i \in A$ and $j \in B \cup C$,
    \item $j < k$ if $j \in B$ and $k \in C$.
\end{itemize}
With respect to this labeling we claim that $\init_{<} J_{G} = \init_{<} J_{H_{1}} + \init_{<} J_{H_{2}}$ where $<$ denotes the lexicographic term order determined by $x_{1} > \cdots x_{n} > y_{1} > \cdots > y_{n}$. By {\cite[Theorem 2.1]{herzog2010binomial}} it suffices to show that if $P$ is an admissible path, i.e. that $P$ is an induced path on the vertices $i = i_{0}, i_{1}, \ldots, i_{\ell} = j$ (where $i_{k}$ is adjacent to $i_{k+1}$) with $i_{k} < i$ or $j < i_{k}$ for all $1 \leq k \leq \ell-1$, then $P \subset H_{1}$ or $P \subset H_{2}$.

First, observe that if $i \in V(A)$, then $j \notin V(C)$. If not, then item \eqref{item:third_item_in_thm_graph_decomp} of Theorem \ref{thm:graph_decomp_implies_conjec} implies that there are no edges of $G$ between the vertices of $A$ and the vertices of $C$. Consequently, $\ell \geq 2$, and there exists $1 \leq k \leq \ell -1$ such that $i_{k} \in B$. But, then our choice of labeling on the vertices would imply that $i < i_{k} < j$ which contradicts $P$ being an admissible path.

Second, observe that if $i \in A$ and $j \in A \cup B$, then $i_{k} \notin C$ for all $1 \leq k \leq \ell-1$. Suppose otherwise, then item \eqref{item:third_item_in_thm_graph_decomp} of Theorem \ref{thm:graph_decomp_implies_conjec} implies that there exists $1 \leq k_{0} < k < k_{s} \leq \ell$ such that the vertices $i_{k_{0}}$ and $i_{k_{s}}$ belong to $B$ and $i_{j} \in C$ for all $k_{0} < j < k_{s}$. Since $B$ induces a connected complete subgraph of $G$ by item \eqref{item:second_item_in_thm_graph_decomp} we have that $\{i_{k_{0}},i_{k_{s}}\}$ is an edge of $G$ which contradicts $P$ being an induced subpath of $G$.

The above two observations show that $\init_{<} J_{G} \subset \init_{<} J_{H_{1}} + \init_{<}J_{H_{2}}$. The inclusion $\supseteq$ follows from the fact that an admissible path of $H_{1}$ remains an admissible path of $G$ since $H_{1}$ is an induced subgraph of $G$. We now have that
\begin{align*}
    \reg(R/J_{G})
    &= \reg(\init_{<} J_{G}) \\
    &\leq \reg(R/\init_{<} J_{H_{1}}) + \reg(R/\init_{<} J_{H_{1}}) \\
    &= \reg(R/J_{H_{1}}) + \reg(R/J_{H_{2}}).
\end{align*}
The first and last equality follow from the fact that binomial edge ideals have squarefree initial ideal with respect to the introduced term order {\cite[Theorem 2.1]{herzog2010binomial}}, and by applying Conca--Varbaro's theorem on preservation of extremal Betti numbers via squarefree Gr\"obner degeneration {\cite[Corollary 2.7]{conca2020square}}. The inequality is Theorem \ref{thm:kalai}.
\eproof

%NRL
\brem
Theorem \ref{thm:graph_decomp_implies_conjec} can be utilized to give a new proof of the computation of $\reg(R/J_{G})$ when $G$ is closed, or when $G$ is bipartite and $J_{G}$ is Cohen--Macaulay. The computation of the Castelnuovo--Mumford regularities for these families of graphs was originally given by Ene--Zarojanu \cite{ene2015regularity} and Jayanthan--Kumar \cite{jayanthan2019regularityCM&BipartiteGraphs}, respectively.
\erem 

\section{Open Questions}

As was mentioned in the introduction it is known for the edge ideal of a graph $G$ that
\begin{align*}
    \reg(I(G)) \leq m(G)+1 \leq \hgt(I(G))+1
\end{align*}
where $m(G)$ denotes the size of a maximal matching of $G$. The notion of size of a maximal path packing of $G$, denoted $r(G)$, provides the analogue of $m(G)$ for binomial edge ideals, see \cite{herzog2022graded}, \cite{laclair2024invariants}. Thus, it is not unreasonable to ask
\bques
For a graph $G$ is
\begin{align*}
    \reg(R/J_{G}) \leq r(G)?
\end{align*}
\eques 
It is known that $im(G) + 1 \leq \reg(I(G))$ where $im(G)$ denotes the size of a maximal induced matching of $G$ \cite{katzman2006BettiNumbers}. The graphs for which $\im(G) = m(G)$ are known as the Cameron--Walker graphs, and admit a combinatorial description \cite{CameronWalker2005Matchings}. An analogue of $im(G)$ for binomial edge ideals is given by $\nu(G)$ in \cite{laclair2024regularity}. It might be interesting to investigate
\bques
Give a combinatorial description for graphs $G$ satisfying $\nu(G) = r(G)$.
\eques

\section*{Acknowledgements}
The author would like to thank Matt Mastroeni for raising a question which led to this paper, and Professors Jason McCullough and Uli Walther for helpful feedback on an earlier version of this paper. LaClair was partially supported by National Science Foundation grant DMS--2100288 and by Simons Foundation Collaboration Grant for Mathematicians \#580839.

%\pagebreak
\bibliographystyle{abbrv}
\bibliography{bibliography}

\end{document}